# Generalized Class of Estimators for Finite Population Mean When Study Variable is Qualitative in Nature

Prayas Sharma, Hemant K.Verma and Rajesh Singh

Department of Statistics, Banaras Hindu University, Varanasi, India.

prayassharma@gmail.com, coolhemant010189@gmail.com ,rsinghstat@gmail.com

**Abstract**

This paper suggests a generalized class of estimators for population mean of the qualitative study variable in simple random sampling using information on an auxiliary variable. Asymptotic expressions of bias and mean square error of the proposed class of estimators have been obtained. Asymptotic optimum estimator has been investigated along with its approximate mean square error. It has been shown that proposed generalized class of estimators are more efficient than all the estimators considered by Singh et al. (2010) in case of qualitative study variable. In addition theoretical findings are supported by an empirical study based on real population to show the superiority of the constructed estimators over others.

1. Introduction

Statisticians are often interested to use auxiliary information in sample surveys at estimation stage in order to improve the precision or accuracy of an estimator of unknown population parameter of interest. In some situations the auxiliary information is not available directly but in the form of an attribute that is auxiliary information is qualitative in nature. When auxiliary information is qualitative in nature then using the point bi-serial correlation between the study variable y and the auxiliary attribute $\phi$ several authors including Naik and Gupta (1996), Jhajj et al. (2006), Singh et al. (2007), Shabbir and Gupta (2007), Singh et al.(2008), Singh et al. (2010), Abd-Elfattah et al. (2010), Singh and Solanki (2012), Sharma et al. (2013a,2013b) proposed improved estimators of population parameters of interest under different situations.

All the authors have implicitly assumed that the study variable Y is quantitative whereas the auxiliary variable is qualitative. But there may be practical situations when study variable itself is qualitative in nature. For example, consider U.S. presidential elections. Assume that there are two political parties, Democratic and Republican. The dependent variable here is the vote choice between two political parties. Suppose we let Y=1, if the vote is for a Democratic candidate and Y=0, if the vote is for republican candidate. Some of the variables used in the vote choice are growth rate of GDP, unemployment and inflation rates, whether the candidate is running for re-election, etc. For the present purposes, the important thing is to note that the study variable is a qualitative variable. One can think several other examples where the study variable is qualitative in nature. Thus, a family either owns a house or it does not, it has disability insurance or it does not, both husband and wife are in the labour force or only one suppose is, etc. In this paper we propose a generalized class estimators in which study variable is qualitative in nature. (see Gujarati and Sangeetha (2007)).

Consider a finite population U= ($U_1$, $U_2$, $U_3$, ...., $U_N$) containing N distinct and identifiable units. Let a sample of size n drawn by simple random sampling without replacement (SRSWOR) from a population U to estimate the population mean of qualitative variable. Let $\phi_i$ and $x_i$ denote the observations on variable $\phi$ and x respectively for $i^{th}$ unit (i=1,2,3...N). $\phi_i = 1$, if $i^{th}$ unit of population possesses attribute $\phi$ and $\phi_i = 0$, otherwise. Further let $A = \sum_{i=1}^{N} \phi_i$ and $a = \sum_{i=1}^{n} \phi_i$, denotes the total number of units in the population and sample possessing attribute $\phi$ respectively, $P = \dfrac{A}{N}$ and $p = \dfrac{a}{n}$, denotes the proportion of units in the population and sample, respectively, possessing attribute $\phi$.

Let us define,

$$e_0 = \dfrac{(p-P)}{P}, \qquad e_1 = \dfrac{(\bar{x}-\bar{X})}{\bar{X}},$$

Such that,

$E(e_i) = 0, (i = 0,1)$

and

$$E(e_0^2) = fC_\phi^2, \qquad E(e_1^2) = fC_x^2, \qquad E(e_1 e_p) = f\rho C_\phi C_x,$$

where,

$$f = \left(\frac{1}{n} - \frac{1}{N}\right), \qquad C_\phi^2 = \frac{S_\phi^2}{\phi^2}, \qquad C_x^2 = \frac{S_x^2}{\overline{X}^2},$$

And $\rho$ is the point bi-serial correlation coefficient between $\phi$ and x.

The remaining portion of this paper is as follows, In section 2 we have considered some existing estimator along with its biases and mean square errors. In section 3, we have suggested a generalized class of estimators along with its members and studied their properties. Section 4, made some comparison of suggested class with other existing estimators. An empirical study is carried out in section 5. We have end the paper with the final conclusion.

## 2. Available estimators in literature when study variable itself an attribute

A ratio-type estimator proposed by Singh et al. (2010) for estimating unknown population mean in case of qualitative study variable is

$$t_S = \left(\frac{P}{\overline{x}}\right)\overline{X} \tag{2.1}$$

The bias and MSE of the estimator $t_S$, to the first order of approximation is given as

$$B(t_S) = f\left(\frac{C_x^2}{2} - \rho C_p C_x\right) \tag{2.2}$$

$$MSE(t_S) = fP^2\left(C_\phi^2 + C_x^2 - 2\rho C_\phi C_x\right) \tag{2.3}$$

Singh et al. (2010) suggested another general class of estimator as,

$$t_{GS} = H(p, u) \tag{2.4}$$

where $u = \dfrac{\overline{x}}{\overline{X}}$ and $H(p,u)$ is a parametric equation of p and u such that

$$H(p,1) = P, \forall P \tag{2.5}$$

and satisfying following regulations:

(i) Whatever be the sample chosen, the point (p,u) assume values in a bounded closed convex subset $R_2$ of the two-dimensional real space containing the point (p,1).

(ii) The function H(p,u) is a continuous and bounded in $R_2$.

(iii) The first and second order partial derivatives of H(p,u) exist and are continuous as well as bounded in $R_2$.

Where,

$$H_1 = \left.\dfrac{\partial H}{\partial u}\right|_{p=P, u=1}, \qquad H_2 = \left.\dfrac{1}{2}\dfrac{\partial^2 H}{\partial u^2}\right|_{p=P, u=1},$$

$$H_3 = \left.\dfrac{1}{2}\dfrac{\partial^2 H}{\partial p \partial u}\right|_{p=P, u=1}, \qquad \text{and} \qquad H_4 = \left.\dfrac{1}{2}\dfrac{\partial^2 H}{\partial p \overline{y}^2}\right|_{p=P, u=1}.$$

The bias and minimum MSE of the estimator $t_b$ are respectively, given by –

$$B(t_{GS}) = f\left(P\rho C_\phi C_x H_3 + C_x^2 H_2 + P^2 C_y^2 H_4\right) \tag{2.6}$$

$$MSE(t_{GS})_{min} = fP^2 C_\phi^2 (1-\rho^2) \tag{2.7}$$

Singh et al. (2010) proposed a new family of estimator for estimating P, as

$$t_{NS} = \left[q_1 P + q_2(\overline{X} - \overline{x})\right]\left[\dfrac{a\overline{X} + b}{a\overline{x} + b}\right]^\alpha \exp\left[\dfrac{(a\overline{X} + b) - (a\overline{x} + b)}{(a\overline{X} + b) + (a\overline{x} + b)}\right]^\beta \tag{2.8}$$

The bias and minimum MSE of the estimator to the first order of approximation, are respectively, given as

$$\text{Bias}(t_{NS}) = P(q-1) + f\left[(q_2\overline{X}B + q_1PA)C_x^2 - q_1P B\rho C_\phi C_x\right] \tag{2.9}$$

$$\text{MSE}(t_{NS})_{min} = \left[P^2 - \frac{\Delta_1\Delta_5^2 + \Delta_3\Delta_4^2 - 2\Delta_2\Delta_4\Delta_5}{\Delta_1\Delta_3 - \Delta_2^2}\right] \tag{2.10}$$

where,

$$M_1 = P^2 f\left(C_p^2 + B^2 C_x^2 - 2B\rho C_p C_x\right), \qquad M_2 = \overline{X}^2 f\left(C_x^2\right),$$

$$M_3 = P^2 f\left(AC_x^2 - 2B\rho C_p C_x\right), \qquad M_4 = P\overline{X}f\left(-BC_x^2 + \rho C_p C_x\right),$$

$$M_5 = \overline{X}Pf\left(-BC_x^2\right)$$

$$q_1^* = \frac{\Delta_1\Delta_4 - \Delta_2\Delta_5}{\Delta_1\Delta_3 - \Delta_2^2} \qquad \text{And} \qquad q_2^* = \frac{\Delta_1\Delta_5 - \Delta_2\Delta_4}{\Delta_1\Delta_3 - \Delta_2^2} \tag{2.11}$$

where,

$$\Delta_1 = (P^2 + M_1 + 2M_3), \Delta_2 = (-M_4 - M_5), \Delta_3 = (M_2), \Delta_4 = (P^2 + M_3), \Delta_5 = (-M_5),$$

## 3. The Suggested Generalised Class of Estimators

We propose a generalized class of estimators for estimating P for a qualitative variable $\phi$, as

$$t_N = \left\{d_1 p\left(\frac{\overline{X}}{\overline{x}}\right)^\alpha \exp\left(\frac{\eta(\overline{X}-\overline{x})}{\eta(\overline{X}+\overline{x})+2\lambda}\right)\right\} + d_2\overline{x} + (1-d_1-d_2)\overline{X} \tag{3.1}$$

where $(d_1, d_2)$ are suitable constants that can be chosen such that MSE of $t_N$ is minimum, $\eta$ and $\lambda$ are either real numbers or the functions of the known parameters of auxiliary variables such as coefficient of variation $C_x$, skewness $\beta_{1(x)}$, kurtosis $\beta_{2(x)}$ and correlation coefficient $\rho$ (see Sharma and Singh (2013c)).

It is to be mentioned that

(i) For $(d_1, d_2) = (1, 0)$, the class of estimator $t_m$ reduces to the class of estimator as

$$t_{NP} = \left\{ p\left(\frac{\overline{X}}{\overline{x}}\right)^\alpha \exp\left(\frac{\eta(\overline{X}-\overline{x})}{\eta(\overline{X}+\overline{x})+2\lambda}\right) \right\} \tag{3.2}$$

(ii) For $(d_1, d_2) = (d_1, 0)$, the class of estimator $t_N$ reduces to the class of estimator as

$$t_{NQ} = \left\{ d_1 p\left(\frac{\overline{X}}{\overline{x}}\right)^\alpha \exp\left(\frac{\eta(\overline{X}-\overline{x})}{\eta(\overline{X}+\overline{x})+2\lambda}\right) \right\} \tag{3.3}$$

Set of new estimators originated from (3.1) choosing the suitable values of $d_1, d_2, \alpha, \eta$ and $\lambda$ are listed in Table 3.1.

**Table 3.1: Set of estimators generated from the class of estimators $t_N$**

| Subset of proposed estimator | $d_1$ | $d_2$ | $\alpha$ | $\eta$ | $\lambda$ |
|---|---|---|---|---|---|
| $t_{N1} = p$ ( usual unbiased estimator) | 1 | 0 | 0 | 0 | 1 |
| $t_{N2} = p\left(\dfrac{\overline{X}}{\overline{x}}\right) = t_S$ (Singh et al. 2010 type) | 1 | 0 | 1 | 0 | 1 |
| $t_{N3} = p\left(\dfrac{\overline{X}}{\overline{x}}\right)^\alpha = \hat{M}_3$ (Srivastava, 1967 type) | 1 | 0 | $\alpha$ | 0 | 1 |
| $t_{N4} = p\left(\dfrac{\overline{x}}{\overline{X}}\right) = M_p$ ( Murthy , 1964 type ) | 1 | 0 | -1 | 0 | 1 |
| $t_{N5} = d_1 p\left(\dfrac{\overline{X}}{\overline{x}}\right)$ (Al and Cingi, 2009 type) | 1 | 0 | 1 | 0 | 1 |
| $t_{N6} = d_1 p\left(\dfrac{\overline{x}}{\overline{X}}\right)$ | $w_1$ | 0 | -1 | 0 | 1 |

| | | | |
|---|---|---|---|
| $t_{N7} = d_1 p$ (Al and Cingi, 2009) | $w_1$ 0 0 0 1 |
| $t_{N8} = w_1 p + w_2 \bar{x} + (1 - w_1 - w_2)\bar{X}$ | $w_1$ $w_2$ 0 0 1 |

Another set of estimators generated from class of estimator $t_{NQ}$ given in (3.3) using suitable values of $\eta$ and $\lambda$ are summarized in table 3.2

### Table 3.2: Set of estimators generated from the estimator $t_{NQ}$

| Subset of proposed estimator | $\alpha$ | $\eta$ | $\lambda$ |
|---|---|---|---|
| $t_{NQ}^{(1)} = \left\{ d_1 p \left(\dfrac{\bar{X}}{\bar{x}}\right) \exp\left(\dfrac{(\bar{X}-\bar{x})}{(\bar{X}+\bar{x})+2}\right) \right\}$ | 1 | 1 | 1 |
| $t_{NQ}^{(2)} = \left\{ d_1 p \left(\dfrac{M_x}{\hat{M}_x}\right) \exp\left(\dfrac{(\bar{X}-\bar{x})}{(\bar{X}+\bar{x})+2\rho_c}\right) \right\}$ | 1 | 1 | $\rho$ |
| $t_{NQ}^{(3)} = \left\{ d_1 p \left(\dfrac{\bar{X}}{\bar{x}}\right) \exp\left(\dfrac{(\bar{X}-\bar{x})}{(\bar{X}+\bar{x})+2\bar{X}}\right) \right\}$ | 1 | 1 | $\bar{X}$ |
| $t_{NQ}^{(4)} = \left\{ d_1 p \left(\dfrac{\bar{X}}{\bar{x}}\right) \exp\left(\dfrac{(\bar{X}-\bar{x})}{(\bar{X}+\bar{x})}\right) \right\}$ | 1 | 1 | 0 |
| $t_{NQ}^{(5)} = \left\{ w_1 p \left(\dfrac{\bar{x}}{\bar{X}}\right) \exp\left(\dfrac{(\bar{X}-\bar{x})}{(\bar{X}+\bar{x})}\right) \right\}$ | -1 | 1 | 1 |
| $t_{NQ}^{(6)} = \left\{ d_1 p \left(\dfrac{\bar{X}}{\bar{x}}\right) \exp\left(\dfrac{\bar{X}(\bar{X}-\bar{x})}{\bar{X}(\bar{X}+\bar{x})+2\rho_c}\right) \right\}$ | 1 | $\bar{X}$ | $\rho$ |
| $t_{NQ}^{(7)} = \left\{ d_1 p \exp\left(\dfrac{\bar{X}(\bar{X}-\bar{x})}{\bar{X}(\bar{X}+\bar{x})+2\rho_c}\right) \right\}$ | 0 | $\bar{X}$ | $\rho$ |
| $t_{NQ}^{(8)} = \left\{ d_1 p \left(\dfrac{\bar{X}}{\bar{x}}\right) \exp\left(\dfrac{\rho_c(\bar{X}-\bar{x})}{\rho_c(\bar{X}+\bar{x})+2\bar{X}}\right) \right\}$ | 1 | $\rho$ | $\bar{X}$ |
| $t_{NQ}^{(9)} = \left\{ d_1 p \left(\dfrac{\bar{x}}{\bar{X}}\right) \exp\left(\dfrac{\rho_c(\bar{X}-\bar{x})}{\rho_c(\bar{X}+\bar{x})+2\bar{X}}\right) \right\}$ | -1 | $\rho$ | $\bar{X}$ |

Expressing (3.1) in terms of e's, we have

$$t_N = d_1 P(1+e_0)(1+e_1)^{-\alpha} \exp\{-ke_1(1+ke_1)^{-1}\} + d_2\overline{X}(1+e_1) + (1-d_1-d_2)\overline{X}$$

where, $k = \dfrac{\eta \overline{X}}{2(\eta \overline{X} + \lambda)}.$ (3.4)

Up to the first order of approximation we have,

$$(t_N - P) = [(d_1 - 1)b + d_1 P\{e_0 - ae_1 + de_1^2 - ae_0 e_1\} + d_2 \overline{X} e_1]$$ (3.5)

where $a = (\alpha + k)$, $b = (P - \overline{X})$ and $d = \left\{\dfrac{3}{2}k^2 + \alpha k + \dfrac{\alpha(\alpha+1)}{2}\right\}$

from equation (3.5), we have

$$(t_N - P)^2 = [(1-2d_1)b^2 + d_1^2\{b^2 + P^2(e_0^2 + a^2 e_1^2 - 2ae_0 e_1)\}$$
$$+ d_2^2 \overline{X}^2 e_1^2 + 2d_1 d_2 P\overline{X}(e_0 e_1 - ae_1^2)]$$

Taking expectations both sides, we get the MSE of the estimator $t_N$ to the first order of approximation as

$$MSE(t_N) = [(1-2w_1)b^2 + d_1^2 M + d_2^2 N + 2d_1 d_2 O]$$ (3.7)

where,

$M = b^2 + P^2 f(C_\phi^2 + a^2 C_x^2 - 2a\rho C_\phi C_x),$

$N = \overline{X}^2 f C_x^2,$

$O = P\overline{X} f(\rho C_\phi - aC_x)C_x.$

The optimum values of $d_1$ and $d_2$ are obtained by minimizing (3.7) and is given by

$$d_1^* = \dfrac{b^2 N}{(MN - O^2)} \quad \text{And} \quad d_2^* = \dfrac{-b^2 O}{(MN - O^2)}$$ (3.8)

Substituting the optimal values of $d_1$ and $d_2$ in equation (3.7) we obtain the minimum MSE of the estimator $t_N$ as

$$MSE_{min}(t_N) = b^2 \left[1 - \frac{b^2 N}{(MN - O^2)}\right] \qquad (3.9)$$

Putting the values of M, N, O and b and simplifying, we get the minimum MSE of estimator $t_N$ as

$$MSE_{min}(t_N) = \left[\frac{P^2(1-R)^2 fC_\phi^2(1-\rho^2)}{[(1-R)^2 + fC_\phi^2(1-\rho^2)]}\right] \qquad (3.10)$$

where $R = \dfrac{\overline{X}}{P}$ and P is defined earlier.

Similarly, the minimum MSE of the class of estimators $t_{NQ}$ is given by

$$MSE_{min}(t_{NQ}) = P^2 \left[\frac{(fC_\phi^2 + a^2\gamma C_x^2 - 2af\rho C_\phi C_x)}{(1 + fC_\phi^2 + a^2 fC_x^2 - 2af\rho C_\phi C_x)}\right] \qquad (3.11)$$

## 4. Efficiency Comparisons

From equations (2.3) and (2.7) we have

$$MSE(t_S) \geq MSE_{min}(t_{GS}) = fP^2(C_\phi^2 + C_x^2 - 2\rho C_\phi C_x) \geq P^2 fC_\phi^2(1-\rho^2)$$

Or $\rho^2 C_\phi^2 + C_x^2 - 2\rho C_\phi C_x \geq 0$ \qquad (4.1)

From equations (2.7) and (3.10) we have

$$MSE_{min}(t_{GS}) \geq MSE_{min}(t_N) = P^2 fC_\phi^2(1-\rho^2) \geq \left[\frac{P^2(1-R)^2 fC_\phi^2(1-\rho^2)}{[(1-R)^2 + fC_\phi^2(1-\rho^2)]}\right]$$

Or $(1-R)^2 + fC_\phi^2(1-\rho^2) \geq (1-R)^2$ \qquad (4.2)

The condition given in (4.2) shows always true.

From equation (2.10) and (3.10) we have

$$MSE_{min}(t_{NS}) \geq MSE_{min}(t_N)$$

If, $\left[P^2 - \dfrac{\Delta_1 \Delta_5^2 + \Delta_3 \Delta_4^2 - 2\Delta_2 \Delta_4 \Delta_5}{\Delta_1 \Delta_3 - \Delta_2^2}\right] \geq \left[\dfrac{P^2(1-R)^2 fC_\phi^2(1-\rho^2)}{[(1-R)^2 + fC_\phi^2(1-\rho^2)]}\right]$ \qquad (4.3)

It follows from (4.1), (4.2) and (4.3), that the proposed class of estimators $t_N$ is better than the ratio estimator $t_S$, general class of estimator $t_{GS}$ and the family of estimators $t_{NS}$, due to Singh et al. (2010) under certain conditions.

**Remark 4.1: Estimator Based on optimum values**

Putting the optimum values of $w_1^*$ and $w_2^*$ in the equation (3.1) we get the optimum estimator as:

$$t'_m = \left\{ d_1^* p \left(\frac{\overline{X}}{\overline{x}}\right)^\alpha \exp\left(\frac{\eta(\overline{X}-\overline{x})}{\eta(\overline{X}+\overline{x})+2\lambda}\right) \right\} + w_2^* \overline{x} + (1 - d_1^* - d_2^*)\overline{X} \qquad (4.4)$$

If the experimenter is not able to specify the value precisely, then it may be desirable to estimate the optimum values from the samples, therefore the values of $w_1^*$ and $w_2^*$ are given as:

$$d_1^* = \frac{\hat{b}^2 \hat{N}}{\hat{M}\hat{N} - \hat{O}^2} \quad \text{and} \quad d_2^* = \frac{\hat{b}^2 \hat{O}}{\hat{M}\hat{N} - \hat{O}^2}$$

where $M = \hat{b}^2 + P^2 \gamma(\hat{C}_y^2 + \hat{a}^2 \hat{C}_x^2 - 2a\hat{\rho}\hat{C}_y \hat{C}_x)$,

$N = \overline{X}\gamma \hat{C}_x^2$, $\hat{\rho}_c = 4(4\hat{p}_{11} - 1)$

$O = P\overline{X}\gamma(\hat{\rho}\hat{C}_\phi - a\hat{C}_x)\hat{C}_x$, $\hat{b} = (P - \overline{X})$, $\hat{a} = (\alpha + \hat{k})$ and $\hat{k} = \frac{\eta \hat{M}_x}{2(\eta \hat{M}_x + \lambda)}$

Here, we have assumed that the population median of auxiliary variable x is known, therefore $\hat{M}_x$ can also be remain as $M_x$.

Expressing (4.4) in terms of e's, we have

$$t'_m = w_1^* P(1+e_0)(1+e_1)^{-\alpha} \exp\left\{-\hat{k}e_1(1+\hat{k}e_1)^{-1}\right\} + w_2^* \overline{X}(1+e_1) + (1 - w_1^* - w_2^*)\overline{X}$$

Proceeding as above, we get the minimum MSE of the estimator $t'_m$ given as:

$$\text{MSE}_{min}(t'_m) = \left[\frac{\hat{M}_y^2 (1-\hat{R})^2 \hat{\gamma}\hat{C}_y^2 (1-\hat{\rho}_c^2)}{[(1-\hat{R})^2 + \hat{\gamma}\hat{C}_y^2(1-\hat{\rho}_c^2)]}\right] \qquad (4.5)$$

## 5. Empirical study

**Data Statistics:** To illustrate the efficiency of proposed generalized class of estimators in the application, we consider the following population data set.

The data used for empirical study has been taken from Gujrati and Sangeetha (2007) -pg, 601. And using raw data we have calculated the following values. Where, y : Home ownership.

x : Income (in thousands of dollars

The values of the required parameters are:

N=40, n=11, $P = 0.525$, $\bar{X} = 14.4$, $C_\phi = 0.963$, $C_x = 0.308$ $\rho = 0.897$, R= 27.42, $\lambda_{12} = -0.118$, $\lambda_{04} = 1.75$, $\lambda_{03} = 0.963$

**Table 5.1: Variances / MSEs/minimum MSEs of different Estimators**

| Estimators | MSE | PRE |
|---|---|---|
| $V(p)$ | 0.061122 | 100.00 |
| $MSE(t_S)$ | 0.32271 | 189.3812 |
| $MSE_{min}(t_{GS})$ | 0.01190 | 511.7912 |
| $MSE_{min}(t_{NS})$ | 0.01171 | 518.9214 |
| $MSE_{min}(t_N)$ | **0.00329** | **1856.8818** |
| $MSE_{min}(t_{N1})$ | 0.01682 | 362.8112 |
| $MSE_{min}(t_{N2})$ | 0.00881 | 687.2571 |
| $MSE_{min}(t_{N3})$ | 0.01191 | 511.7912 |
| $MSE_{min}(t_{N4})$ | 0.02801 | 216.3089 |
| $MSE_{min}(t_{N5})$ | 0.00881 | 687.2763 |
| $MSE_{min}(t_{N6})$ | 0.02821 | 216.3019 |
| $MSE_{min}(t_{N7})$ | 0.01681 | 362.8229 |
| $MSE_{min}(t_{N8})$ | **0.00329** | **1856.8818** |
| $MSE_{min}(t_{NQ}^1)$ | 0.00636 | 960.8345 |
| $MSE_{min}(t_{NQ}^2)$ | 0.00631 | 963.0277 |
| $MSE_{min}(t_{NQ}^3)$ | 0.00744 | 820.9345 |
| $MSE_{min}(t_{NQ}^4)$ | 0.00621 | 983.6847 |
| $MSE_{min}(t_{NQ}^5)$ | 0.02211 | 276.3287 |
| $MSE_{min}(t_{NQ}^6)$ | 0.00622 | 982.1553 |
| $MSE_{min}(t_{NQ}^7)$ | 0.01245 | 490.7537 |
| $MSE_{min}(t_{NQ}^8)$ | 0.00151 | 812.9560 |
| $MSE_{min}(t_{NQ}^9)$ | 0.02521 | 242.0966 |

Table 5.1 exhibits variance, mean square errors and percent relative efficiencies of the existing estimators $p$, $t_S$, $t_{GS}$, $t_{NS}$ and proposed generalised class of estimators along with its different members. Analysing table 5.1 we conclude that the estimators based on auxiliary information are more efficient than the one which does not use the auxiliary information as $p$. The members $t_{NQ}^1$, $t_{NQ}^2$, $t_{NQ}^4$ and $t_{NQ}^6$ of the class of estimators $t_{NQ}$, obtained from generalized class of estimators $t_N$, are almost equally efficient but more than the usual unbiased estimator $p$, usual ratio estimator $t_S$, class of estimators $t_{GS}$ and family of estimators $t_{NS}$. Among the proposed class of estimators $t_N$ and $t_{NQ}^j$ (j=1,2,...9) the performance of the estimator $t_N$, which is equal efficient to the estimator $t_{N8}$, is best in the sense of having the least MSE (maximum PRE) followed by the estimator $t_{NQ}^4$ and $t_{NQ}^6$ which utilize the information of population mean $\overline{X}$ and correlation coefficient $\rho$.

**Conclusion**

In this article we have suggested a generalized class of estimators for the population mean when study variable is qualitative in nature using auxiliary information in simple random sampling without replacement (SRSWOR). In addition, some known estimators of population mean when study variable is an attribute such as usual unbiased estimator p, ratio estimator due to Singh et al.(2010) , Srivastava (1967) type estimator , murthy(1964) type estimator, Al and Chingi (2009) and Singh and Solanki (2013) type estimators are found to be members of the proposed generalized class of estimators. Some new members are also generated from the proposed generalized class of estimators using auxiliary information. We have determined the biases and mean square errors of the proposed class of estimators up to the first order of approximation. The proposed generalized class of estimators are advantageous in the sense that the properties of the estimators, which are members of the proposed class of estimators, can be easily obtained from the properties of the proposed generalized class. Thus the study unifies properties of several estimators for population mean when study variable is qualitative in nature. In theoretical and empirical efficiency comparisons, it has been shown that the proposed generalized class of estimators $t_N$ are more efficient than the estimators considered here.